\documentclass[praybox]{svmult}

\usepackage{mathptmx}       
\usepackage{helvet}         
\usepackage{courier}        

\usepackage{makeidx}         
\usepackage{graphicx}        
\usepackage{multicol}        
\usepackage[bottom]{footmisc}
\usepackage{amsmath}
\usepackage{amssymb}

\usepackage{url}

\usepackage{tikz}

\newcommand{\BBM}{\begin{bmatrix}}
\newcommand{\EBM}{\end{bmatrix}}
\newcommand{\BEAS}{\begin{eqnarray*}}
\newcommand{\EEAS}{\end{eqnarray*}}
\newcommand{\BEA}{\begin{eqnarray}}
\newcommand{\EEA}{\end{eqnarray}}
\newcommand{\BEQ}{\begin{equation}}
\newcommand{\EEQ}{\end{equation}}
\newcommand{\BIT}{\begin{itemize}}
\newcommand{\EIT}{\end{itemize}}

\newcommand{\reals}{{\mbox{\bf R}}}
\newcommand{\naturals}{{\mbox{\bf N}}}

\newcommand{\diag}{\mathop{\bf diag}}

\DeclareMathOperator*{\Ne}{Ne}
\DeclareMathOperator*{\argmin}{argmin}
\DeclareMathOperator*{\children}{ch}
\DeclareMathOperator*{\parent}{par}

\makeindex    

\begin{document}

\title*{Exploiting Chordality in Optimization Algorithms for Model 
Predictive Control}

\author{Anders Hansson and 
Sina Khoshfetrat Pakazad}

\institute{Anders Hansson 
\at Division of Automatic Control, Linkoping University, Sweden, \email{anders.g.hansson@liu.se}
\and Sina Khoshfetrat Pakazad 
\at C3 IoT, Redwood city, California, USA, 
\email{sina.pakazad@c3iot.com}}

\maketitle
\abstract{In this chapter we show that chordal structure can be used
to devise efficient optimization methods for many common model predictive
control problems. The chordal structure is used both for 
computing search directions efficiently as well as for distributing all the
other computations in an interior-point method for solving the problem. 
The chordal structure can stem both from the sequential nature of the problem as well as from distributed formulations of the problem related to scenario trees
or other formulations. The framework enables efficient parallel computations. 
}

\section{Introduction}
Model Predictive Control (MPC) is an important class of controllers that 
are being employed more and more in industry, 
\cite{qin+bad03}. It has its root
going back to \cite{cut+ram79}. The success is mainly because it can handle 
constraints on control signals and/or states in a systematic way. In the
early years its applicability was limited to slow processes, since 
an optimization problem has to be solved at each sampling instant. 
Tremendous amount of research has been spent on overcoming this limitation. 
One avenue has been what is called explicit MPC, \cite{Alessio2009}, 
where the optimization
problem is solved parametrically off-line. Another avenue has been to 
exploit the inherent structure of the optimization problems stemming from MPC, 
\cite{gla+jon84,wri93,ste94,arn+put94,wri96,rao+wri+raw97,gop+bie98,han98,van+boy+nou01,jor04,ake+han04,Diehl2009,axe+van+han10,wan+boy10,jer+ker+con12,dom+12,fri15,kli17,nie17}. Typically
this has been to use Riccati recursions to efficiently compute search directions for 
Interior Point (IP) methods or actives set methods to solve the optimization problem.
In this paper we will argue that the important structures that have been 
exploited can all be summarized as {\it chordal structure}. Because of this
the same structure exploiting software can be used to speed up all computations
for MPC. This is irrespective of what MPC formulation is considered and 
irrespective of what type of optimization algorithm is used. 
We assume that the reader is familiar with the receding horizon strategy
of MPC and we will only discuss the associated constrained finite-time
optimal control problems. We will from now on refer to 
the associated problem as the MPC problem.
We will mostly assume quadratic cost and linear
dynamics and inequality constraints. Even if not all problems fall into
this category, problems with quadratic objective and linear constraints
are often solved as subproblems in solvers.

The remaining part of the paper is organized as follows.
We will in Section~2 discuss how chordal sparsity arises and how it can be 
utilized in general convex optimization problems to obtain 
computations distributed over a so called clique tree. The presentation is based
on \cite{khoshfetrat2016distributed}.
In Section~3 we will then discuss the classical formulation
of MPC, and how the problem can be solved using an IP
method. Specifically we will discuss how the equations for the search
directions can be distributed over the clique tree. The well-known backward
dynamic programming solution will be derived as a special case. We will
see that we can also do forward dynamic programming, combinations of
forward and backward dynamic programming, and even dynamic programming
in parallel. In Section~4 we will discuss regularized MPC. In Section~5
we will discuss stochastic MPC. In Section~6 we will discuss
distributed MPC, and finally in Section~7 we will give some conclusions, 
discuss generalizations of our results 
and directions for future research. 
\subsection*{Notation}
We denote with $\reals$ the set of real numbers, with $\reals^n$ the set of 
$n$-dimensional real-valued vectors and with $\reals^{m\times n}$ the set of 
real-valued matrices with $m$ rows and $n$ columns. 
We denote by $\naturals$ the set of natural
numbers and by $\naturals_n$ the subset  $\lbrace 1,2,\ldots,n\rbrace$ of $\naturals$.
For a vector $x\in\reals^n$ the matrix $X=\diag(x)$ is a diagonal matrix with the
components of $x$ on the diagonal. For two matrices $A$ and $B$ the matrix $A\oplus B$ is a block-diagonal matrix with 
$A$ as the 1,1-block and $B$ as the 2,2-block. For a symmetric matrix $A$ 
the notation $A(\succeq)\succ 0$ is equivalent to $A$ being positive (semi)-definite. 
\section{Chordal Sparsity and Convex Optimization}
Consider the following convex optimization problem
\begin{align}\label{eqPaperIII:CP}
\min_{x} \quad F_1(x) + \dots + F_N(x),
\end{align}
where $F_i \ : \ \reals^{n}\rightarrow \reals$ for all $i = 1, \dots, N$. We assume that each function $F_i$ is only dependent on a small subset of elements of $x$. Let us denote the ordered set of these indexes by $J_i \subseteq \naturals_n$. We also denote the ordered set of indexes of functions that depend on $x_i$ with $\mathcal I_i = \{ k \ | \ i \in J_k\} \subseteq \naturals_N$. We can then rewrite the problem in~\eqref{eqPaperIII:CP}, as
\begin{align}\label{eqPaperIII:CPS}
\min_{x} & \quad   \bar F_1(E_{J_1}x) + \dots + \bar F_N(E_{J_N}x),
\end{align}
where $E_{J_i}$ is a $0$--$1$ matrix that is obtained from an identity matrix of order $n$ by deleting the rows indexed by $\naturals_n \setminus J_i$.  The functions $\bar F_i \ : \ \reals^{|J_i|} \rightarrow \reals$ are lower dimensional descriptions of $F_i$s such that $F_i(x) = \bar F_i(E_{J_i}x)$ for all $x \in \reals^n$ and $i \in\naturals_N$. For instance consider the following optimization problem
\begin{align}
\min_x \quad F_1(x) + F_2(x) +  F_3(x) + F_4(x) + F_5(x) + F_6(x),
\end{align}
and let us assume that $x \in \reals^8$, $J_1 = \{ 1, 3 \}$, $J_2 = \{ 1, 2, 4 \}$, $J_3 = \{ 4, 5 \}$, $J_4 = \{ 3, 4 \}$, $J_5 = \{ 3, 6, 7 \}$ and $J_6 = \{ 3, 8 \}$. We then have $\mathcal I_1 = \{ 1, 2 \}$, $\mathcal I_2 = \{ 2 \}$, $\mathcal I_3 = \{ 1, 4, 5, 6 \}$, $\mathcal I_4 = \{ 2, 3, 4 \}$, $\mathcal I_5 = \{ 3 \}$, $\mathcal I_6 = \{ 5 \}$, $\mathcal I_7 = \{ 5 \}$  and $\mathcal I_8 = \{ 6 \}$.  This problem can then be written in the same format as in \eqref{eqPaperIII:CPS} as
\begin{equation}
\min_x \quad \bar F_1(x_1,x_3) + \bar F_2(x_1, x_2, x_4) +\\ \bar F_3(x_4, x_5) + \bar F_4(x_3, x_4) + \bar F_5(x_3, x_6, x_7) + \bar F_6(x_3, x_8).\label{eqPaperIII:example}
\end{equation}
The formulation of coupled problems as in \eqref{eqPaperIII:CPS} enables us to get a more clear picture of the coupling in the problem. Next we describe how the coupling structure in~\eqref{eqPaperIII:CP} can be expressed graphically using undirected graphs.
\subsection{Sparsity Graph}
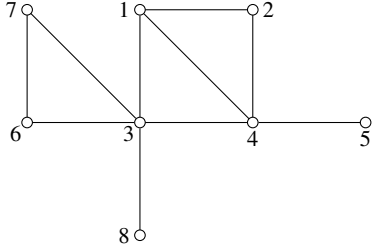
\begin{figure}
		
		\centering
		\begin{tikzpicture}
		\tikzstyle{every node}=[draw,circle,fill=white,minimum size=4pt,
		inner sep=0pt]
		
		\draw (0,0) node (1) [label=left:$1$] {}
		-- ++(0:1.5cm) node (2) [label=right:$2$] {}
		-- ++(-90:1.5cm) node (4) [label=below:$4$] {}
		-- (1);

        \draw (1) 
        -- ++(-90:1.5cm) node (3) [label=below left:$3$] {}
        -- (4);
        
        \draw (3) 
        -- ++(180:1.5cm) node (6) [label=below left:$6$] {}
        -- ++(90:1.5cm) node (7) [label=left:$7$] {}
        -- (3);
        
        \draw (3) 
        -- ++(-90:1.5cm) node (8) [label=left:$8$] {}
        ;
        
        \draw (4) 
        -- ++(0:1.5cm) node (5) [label=below:$5$] {}
        ;
        
		\end{tikzpicture}

\caption{\small The sparsity graph for the problem in~\eqref{eqPaperIII:example}.\normalsize }
\label{figPaperIII:SG}		
		
	\end{figure}

A graph $G$ is specified by its vertex and edge sets $V$ and $\mathcal E$, respectively. A graph that sheds light on the coupling structure of the problem is the so-called \emph{sparsity graph}, $G_s$, of the problem. This graph is undirected with vertex set $V_s = \naturals_n$ and edge set $\mathcal{E}_s$ with $(i,j) \in \mathcal E_s$ if and only if $\mathcal I_i \cap \mathcal I_j \neq \emptyset$. Let us now reconsider the example in~\eqref{eqPaperIII:example}. The sparsity graph for this problem is illustrated in Figure~\ref{figPaperIII:SG}. 

As we will see later graph representations of the coupling structure in problems play an important role in designing distributed algorithms for solving coupled problems and gaining insight regarding their distributed implementations. Specifically, chordal graphs and their characteristics play a major role in the design of our proposed algorithm. This is the topic of the next subsection.
\subsection{Chordal Graphs and Clique Trees}\label{secPaperIII:chordal}
A graph $G(V,\mathcal E)$ with vertex set $V$ and edge set $\mathcal E$ is chordal if every of its cycles of length at least four has a chord, where a chord is an edge between two non-consecutive vertices in a cycle, \cite[Ch. 4]{gol:04}. A clique of $G$ is a \emph{maximal} subset of $V$ that induces a complete subgraph on $G$. Consequently, no clique of $G$ is entirely contained in any other clique, \cite{blp:94}. Let us denote the set of cliques of $G$ as $\mathbf C_G = \{ C_1, \dots, C_q \}$. There exists a tree defined on $\mathbf C_G $ such that for every $C_i, C_j \in \mathbf C_G $ with $i \neq j$, $C_i \cap C_j$ is contained in all the cliques in the path connecting the two cliques in the tree. This property is called the \emph{clique intersection property}, and trees with this property are referred to as clique trees. For instance the graph in Figure \ref{figPaperIII:SG} is chordal and has five cliques, namely $C_1 = \{ 1, 2, 4 \}$, $C_2 = \{ 1, 3, 4 \}$, $C_3 = \{ 4, 5 \}$, $C_4 = \{ 3, 6, 7 \}$ and $C_5 = \{ 3, 8 \}$. A clique tree over these cliques is given in Figure \ref{figPaperIII:SC}. This tree then satisfies the clique intersection property, e.g., notice that $C_2 \cap C_3 = \{ 4 \}$ and the only clique in the path between $C_2$ and $C_3$, that is $C_1$, also includes $\{ 4 \}$.

Chordal graphs and their corresponding clique trees play a central role in our distributed algorithm. For chordal graphs there are efficient methods for computing cliques and clique trees. Sparsity graphs do not have to be chordal. However, there are simple heuristic methods, \cite{cor:01,kol:09}, to compute a chordal embedding of such graphs, where a chordal embedding of a graph $G(V,\mathcal E)$ is a chordal graph with the same vertex set and an edge set $\mathcal E_e$ such that $\mathcal E \subseteq \mathcal E_e$. For the MPC problems
we consider we will derive chordal embeddings manually when required. 
We will now discuss distributed optimization using message-passing.
\subsection{Distributed Optimization Using Message-passing}
\begin{figure}
		
		\centering
		\begin{tikzpicture}
		
		
		\tikzstyle{every node}=[draw,circle,fill=white,minimum size=2pt,
		inner sep=2pt]
		
		\draw (0,0) node (1) [label=right:{$C_1=\{1,2,4\}$}] {1}
		-- ++(225:1.5cm) node (2) [label=left:{$C_2=\{1,3,4\}$}] {2}
		-- ++(225:1.5cm) node (4) [label=left:{$C_4=\{3,6,7\}$}] {4} 
		;
		
        \draw (2) 
		-- ++(-45:1.5cm) node (5) [label=right:{$C_5=\{3,8\}$}] {5}
        ;
		
        \draw (1) 
		-- ++(-45:1.5cm) node (3) [label=right:{$C_3=\{4,5\}$}] {3}  
		;

		\end{tikzpicture}
\caption{\small Clique tree for the sparsity graph of the problem  in \eqref{eqPaperIII:example}.\normalsize }
\label{figPaperIII:SC}		
	\end{figure}
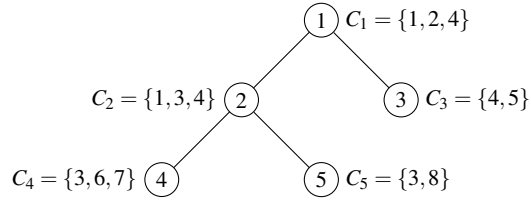

Consider the optimization problem in \eqref{eqPaperIII:CP}. Let $G_s(V_s,\mathcal E_s)$ denote the chordal sparsity graph for this problem and let $\mathbf C_s = \{ C_1, \dots, C_q  \}$ and $T(V_t, \mathcal E_t)$ be its set of cliques and a corresponding clique tree, respectively. It is possible to devise a distributed algorithm for solving this problem that utilizes the clique tree $T$ as its computational graph. This means that the nodes $V_t = \naturals_q$ act as computational agents and collaborate with their neighbors that are defined by the edge set $\mathcal E_t$ of the tree. For example, the sparsity graph for the problem in \eqref{eqPaperIII:example} has five cliques and a clique tree over these cliques is illustrated in Figure~\ref{figPaperIII:SC}. This means the problem can be solved distributedly using a network of five computational agents, each of which needs to collaborate with its neighbors as defined by the edges of the tree, e.g., Agent 2 needs to collaborate with agents $1, 4, 5$.

In order to specify the messages exchanged among these agents, we first assign different terms of the objective function in \eqref{eqPaperIII:CP} to each agent. A valid assignment is that $F_i$ can only be assigned to agent $j$ if $J_i \subseteq C_j$. We denote the ordered set of indices of terms of the objective function assigned to agent $j$ by $\phi_j$. For instance, for the problem in \eqref{eqPaperIII:example}, assigning $\bar F_1$ and $\bar F_4$ to Agent 2 would be a valid assignment since $J_1, J_4 \subseteq C_2$ and hence $\phi_2 = \{ 1, 4 \}$. Notice that the assignments are not unique and for instance there can exist agents $j$ and $k$ with $j \neq k$ so that $J_i \subseteq C_j$ and $J_i \subseteq C_k$ making assigning $F_i$ to agents $j$ or $k$ both valid. It can be shown that for every term of the objective function there will always exist an agent that it can be assigned to.

We will now express the messages that are exchanged among neighboring agents. Particularly, let $i$ and $j$ be two neighboring agents, then the message sent from agent $i$ to agent $j$, $m_{ij}$, is given by
\begin{align}\label{eqPaperIII:mij}
m_{ij}(x_{_{S_{ij}}}) = \min_{x_{_{C_i \setminus S_{ij}}}} \left\{  \sum_{k \in \phi_i} \bar F_k(x_{_{J_k}}) +  \sum_{k \in \Ne(i)\setminus \{j\}} m_{ki}(x_{_{S_{ik}}}) \right\},
\end{align}
where $S_{ij} = C_i \cap C_j$ is the so-called separator set of agents $i$ and $j$. As a result, for agent $i$ to be able to send the correct message to agent $j$ it needs to wait until it has received all the messages from its neighboring agents other than $j$. Hence, the information required for computing a message also sets the communication protocol for this algorithm. Specifically, it sets the ordering of agents in the message-passing procedure in the algorithm, where messages can only be initiated  from the leaves of the clique tree and upwards to the root of the tree, which is referred to as an upward pass through the tree. For instance, for the problem in \eqref{eqPaperIII:example} and as can be seen in Figure~\ref{figPaperIII:SC}, $\Ne(2) = \{ 1, 4, 5 \}$. Then the message to be sent from Agent 2 to Agent 1 can be written as
\begin{align}\label{eqPaperIII:mijexample}
m_{21}(x_1, x_4) = \min_{x_3} \left\{  \bar F_1(x_1, x_3) + \bar F_4(x_3, x_4) +  m_{42}(x_3) + m_{52}(x_3) \right\}.
\end{align}
which can only be computed if Agent 2 has received the messages from agents 4 and 5.

The message, $m_{ij}$, that every agent $j$ receives from a neighboring agent $i$ in fact summarizes all the necessary information that agent $j$ needs from all the agents on the $i$-side of the edge $(i, j)$. 
With this description of messages and at the end of an upward-pass through the clique tree, the agent at the root of the tree, indexed $r$, will have received messages from all its neighbors. Consequently, it will have all the necessary information to compute its optimal solution by solving the following optimization problem
\begin{align}\label{eqPaperIII:RLocalProblem}
x^\ast_{_{C_r}} = \argmin_{x_{_{C_r}}} \left\{  \sum_{k \in \phi_r} \bar F_k(x_{_{J_k}})  + \sum_{k \in \Ne(r)} m_{kr}(x_{_{S_{rk}}}) \right\}.
\end{align}
Then the root sends the computed optimal solution $\left(x^\ast_{_{S_{rj}}}\right)^r$ to its children, i.e., to all agents $j \in  \children(r)$. Here $\left(x^\ast_{_{S_{rj}}}\right)^r$ denotes the optimal solution computed by agent $r$. Then all these agents, similar to the agent at the root, will then have received messages from all their neighbors and can compute their corresponding optimal solution as
\begin{align}\label{eqPaperIII:LocalProblemi}
x^\ast_{_{C_i}} = \argmin_{x_{_{C_i}}} \left\{  \sum_{k \in \phi_i} \bar F_k(x_{_{J_k}})  + \sum_{k \in \Ne(i)\setminus r} m_{ki}(x_{_{S_{ik}}}) +  \right\}.
\end{align}
The same procedure is executed downward through the tree until we reach the leaves, where each agent $i$, having received the computed optimal solution by its parent, i.e., $\left( x_{_{S_{\parent(i)i}}}^\ast \right)^{\parent(i)}$, computes its optimal solution by
\begin{align}\label{eqPaperIII:LocalProblempar}
x^\ast_{_{C_i}} = \argmin_{x_{_{C_i}}} \left\{  \sum_{k \in \phi_i} \bar F_k(x_{_{J_k}})   + \sum_{k \in \Ne(i)\setminus \parent(i)} m_{ki}(x_{_{S_{ik}}}) \right\}.
\end{align}
where $\parent(i)$ denotes the index for the parent of agent $i$.

Notice that in case the optimal solution of \eqref{eqPaperIII:CP} is not unique, then we need to modify the algorithm with 
regularization terms, see \cite{khoshfetrat2016distributed}.

So far we have provided a distributed algorithm to compute optimal solution  for convex optimization problems in the form \eqref{eqPaperIII:CP}. However, this algorithm relies on the fact that we are able to eliminate variables and compute the optimal objective value as a function of the remaining ones in closed form. This capability is essential, particularly for computing the exchanged messages among agents and in turn seemingly limits the scope of problems that can be solved using this algorithm. It turns out that 
the described algorithm can be incorporated within a primal-dual interior-point method to solve general convex optimization problems, distributedly, \cite{khoshfetrat2016distributed}. Then the message passing
algorithm is used to solve
the quadratic subproblems that approximate the overall problem at
each and every iterate in order to compute search directions. Hence the messages will be quadratic functions that are easy to represent. 
\subsection{Interior-Point Methods}
All the MPC problems that we encountere in this chapter
are special cases of Quadratic Programs (QPs). We will now discuss
how such problems can be solved using IP methods, \cite{wri97}. Consider the QP
\begin{align}
\min_{z} &\frac{1}{2}z^T\mathcal Q z+q^Tz\label{QP:objective}\\
{\rm s.t.}\;& \mathcal A z = b\label{QP:eq}\\
&\mathcal D z\leq e\label{QP:ineq}
\end{align}
where $\mathcal Q\succeq 0$, i.e. positive semidefinite, where $\mathcal A$ has full row rank, and where the matrices and vectors are of compatible dimensions. The inequality in (\ref{QP:ineq})
is component-wise inequality. 
The Karush-Kuhn-Tucker (KKT) optimality conditions for this problem is 
\begin{equation}
\begin{bmatrix}\mathcal Q&\mathcal A^T&\mathcal D^T&\\
\mathcal A & & &\\
\mathcal D& & &I\\
& & & M\end{bmatrix}\begin{bmatrix}z\\\lambda\\\mu\\s\end{bmatrix}=
\begin{bmatrix}-q\\b\\e\\0\end{bmatrix}\label{eqn:KKT_classical}
\end{equation}
and $(\mu,s)\geq 0$, where $M=\diag(\mu)$. Blank entries in a matrix
are the same as zero entries. Above $\lambda$ and $\mu$
are the Lagrange multipliers for the equality and inequality constraints, respectively. The vector $s$ is the slack variable
for the inequality constraints. In IP methods one linearizes
the above equations to obtain equations for search directions:
\begin{equation}
\begin{bmatrix}\mathcal Q&\mathcal A^T&\mathcal D^T&\\
\mathcal A & & &\\
\mathcal D& & &I\\
& & S& M\end{bmatrix}\begin{bmatrix}\Delta z\\\Delta\lambda\\\Delta\mu\\\Delta s\end{bmatrix}=
\begin{bmatrix}r_z\\r_\lambda\\r_\mu\\r_s\end{bmatrix}
\end{equation}
where $S=\diag(s)$, and where $r=(r_z,r_\lambda,r_\mu,r_s)$ is some
residual vector that depends on what IP method is used. The
quantities $r$, $S$ and $M$ depend on the value of the current 
iterate in the IP method. From the last two 
rows above we have $\Delta s=r_\mu-\mathcal D\Delta z$ and $\Delta\mu
=S^{-1}(r_s-M\Delta s)$. After substitution of these expressions into the first
two rows we obtain
\begin{equation}
\begin{bmatrix}\mathcal Q+\mathcal D^TS^{-1}M\mathcal D&\mathcal A^T\\
\mathcal A & \end{bmatrix}\begin{bmatrix}\Delta z\\\Delta\lambda\end{bmatrix}=
\begin{bmatrix}r_z-\mathcal D^TS^{-1}(r_s-Mr_\mu)\\r_\lambda\end{bmatrix}\label{eqn:KKT_classical_reduced}
\end{equation}
We notice that the search directions are 
obtained by solving an indefinite symmetric
linear system of eqations. The indefinite 
matrix is reffered to as the KKT matrix, and it is invertible if and only
if 
\begin{equation}
\mathcal Q_s=\mathcal Q+\mathcal D^TS^{-1}\mathcal D\label{eqn:Qs}
\end{equation}
is positive definite on the nullspace of $\mathcal A$. 
Notice that the KKT matrix for the search directions can be interpreted as the optimality conditions of a QP with only equality constraints, where the quadratic weight is modified such that it is larger the closer the iterates are to the boundary of the
constraints. In case this QP is loosely coupled with chordal 
structure message passing over a clique tree can be used to compute
the search directions in a distributed way as described above. 
The key to this will be to solve parametric QPs, which is the next
topic. Before finishing this section we remark that for active set methods
similar QPs also have to be solved. There will however be additional 
equality constraints depending on what constraints are active at the current
iterate. 
\subsection{Parametric QPs}
Consider the quadratic optimization problem 
\begin{align}
\min_{z} &\frac{1}{2}z^TMz+m^Tz\\
{\rm s.t.}\;& Cz=d
\end{align}
with $C$ full row rank and $M\succeq 0$. The KKT conditions for the optimal solution are:
$$\begin{bmatrix}M&C^T\\C&\end{bmatrix}
\begin{bmatrix}z\\\lambda\end{bmatrix}=
\begin{bmatrix}-m\\d\end{bmatrix}$$
These equations have a unique solution if and only if $M+C^TC\succ 0$. 
Now consider the partioning of the above problem defined by
$$M = \begin{bmatrix}Q&S\\S^T&R\end{bmatrix};\quad
C = \begin{bmatrix}A& B\\ &D\end{bmatrix};\quad 
d =\begin{bmatrix}e\\f\end{bmatrix};\quad
m =\begin{bmatrix}q\\r\end{bmatrix};\quad
z = \begin{bmatrix}x\\y\end{bmatrix}$$
with $A$ full row rank. We then assume that the variables related to 
a specific leaf are $x$, and we want to solve the problem associated
with the leaf, i.e. 
\begin{align}
\min_{x} &\frac{1}{2}\begin{bmatrix}x\\y\end{bmatrix}^T
\begin{bmatrix}Q&S\\S^T& \end{bmatrix}
\begin{bmatrix}x\\y\end{bmatrix}+q^Tx\\
{\rm s.t.}\;& Ax+By=e
\end{align}
with respect to all $y$. The KKT conditions for this problem
are
$$\begin{bmatrix}Q&A^T\\A&\end{bmatrix}
\begin{bmatrix}x\\\mu\end{bmatrix}=
\begin{bmatrix}-q-Sy\\e-By\end{bmatrix}$$
Notice that the solution $x$ will be affine in $y$, and hence when 
it is substituted back into the objective function we obtain
a quadratic message in $y$. 
The 1,1-block of $M+C^TC$ is $Q+A^TA$, which by the Schur complement
formula is positive definite. Hence the leaf problem has a unique
solution. If we then substitute the solution of the leaf into the
overall problem, we will have a unique solution also for this
problem, since the overall problem has a unique solution. Beceause
of this, every leaf in the message passing algorithm will have a 
problem with unique solutions assuming that the overall problem 
has a unique solution. This also goes for all nodes, since they will 
become leaves as other leaves are pruned away. 

Notice that it is always possible to make
sure that the matrix $A$ has full rank for a leaf
by pre-processing of the inequality constraints. In case $A$ does
not have full row rank, perform a rank-revealing factorization
such that the constraints can be written
$$\begin{bmatrix}\bar A_1\\0\end{bmatrix}x+
\begin{bmatrix}\bar B_1\\\bar B_2\end{bmatrix}y=
\begin{bmatrix}\bar e_1\\\bar e_2\end{bmatrix}$$
and append the constraint $\bar B_2y=\bar e_2$ to belong to 
$$Dy=f$$
This can be done recursivley over the clique tree so that the 
parametric QPs for each node satisfies the rank condition, 
\cite{khoshfetrat2016distributed}.

We will now see how chordal sparsity and distributed computations can
be used to solve optimization problems arising in MPC efficiently. 
\section{Classical MPC}
A classical MPC problem can be cast in the form
\begin{align}
\min_{u} &\frac{1}{2}\sum_{k=0}^{N-1}\begin{bmatrix}x_k\\u_k\end{bmatrix}^TQ
\begin{bmatrix}x_k\\u_k\end{bmatrix}+\frac{1}{2}x_N^TSx_N
\label{problem:classical_objective}\\
{\rm s.t.}\;& x_{k+1}=Ax_k+Bu_k+v_k,\quad x_0=\bar x
\label{problem:classical_dynamics}\\
&Cx_k+Du_k\leq e_k\label{problem:classical_constraints}
\end{align}
with $A$, $B$, $C$, $D$, $Q$, $S$, $\bar x$, $e_k$, and $v_k$ given, and where 
$u=(u_0,u_1,\ldots,u_{N-1})$ are the optimzation variables. The dimensions 
of the control signal $u_k$ and the state vector $x_k$ are $m$ and $n$, 
respectively. The number of inequality constraints are $q$ for each
time index $k$. The dimensions of all other quantities are defined to
be consistent with this. 
We assume that $Q\succeq 0$ and that $S\succeq 0$. This is a convex
quadratic optimization problem. When the inequality constraints are not
present it is a classical Linear Quadratic (LQ) control problem. It 
is of course possible to extend the problem formulation to time-varying
dynamics, inequality constraints and weights. Also the extention to a linear
term in the objective function is straight forward. 
\subsection{Quadratic Program}
The classical formulation in 
(\ref{problem:classical_objective}--\ref{problem:classical_constraints}) 
is equivalent to (\ref{QP:objective}--\ref{QP:ineq})
with $q=0$, 
\begin{align*}
z&=(x_0,u_0,x_1,u_1\ldots,x_{N-1},u_{N-1},x_N)\\
\lambda&=(\lambda_0,\lambda_1,\ldots,\lambda_N)\\
b&=(\bar x,v_0,v_1,\ldots,v_{N-1})\\
e&=(e_0,e_1,\ldots,e_{N-1})
\end{align*}
and 
\begin{align*}
\mathcal A&=
\begin{bmatrix}
 I &   &   &   &  &      &   &   &  \\
 -A&-B &I  &   &  &      &   &   &  \\
   &   &-A &-B &I &      &   &   &  \\
   &   &   &   &  &\ddots&   &   &  \\
   &   &   &   &  &      &-A &-B &I \\
\end{bmatrix}\\
\mathcal D&=\begin{bmatrix}C&D\end{bmatrix}\oplus
\begin{bmatrix}C&D\end{bmatrix}\oplus\ldots\oplus
\begin{bmatrix}C&D\end{bmatrix}\\
\mathcal Q &= Q\oplus Q\oplus\cdots\oplus Q\oplus S
\end{align*}
We see that the data matrices are banded. Hence, sparse linear system 
solvers could be used  when solving the KKT equations for 
search directions in an IP method, but we will see that the structure within the bands
can be further utilized. Also we notice that the matrix $\mathcal Q_s$ in (\ref{eqn:Qs}) has
the same structure as $\mathcal Q$. Therefore the KKT matrix for the search directions can be interpreted as the optimality conditions of an unconstrained LQ control problem for the search directions, where the weights are modified such that they are larger the closer the iterates are to the boundary of the
constraints. 
Since inequality constraints not affect the structure of the KKT matrix we will 
from now on not consider them when we discuss
the different MPC problems. 

The classical formulation without any constraints is usually solved using
a backward Riccati recursion. We will see how this can be obtained from the 
general techniques presented above. This is the same derivation that is
usually done using backward dynamic programming. We will also investigate
forward dynamic programming, and finally we will se how one can obtain 
parallel computations. 
\subsection{Backward Dynamic Programming}
\begin{figure}
		
		\centering
		\begin{tikzpicture}
		\tikzstyle{every node}=[draw,circle,fill=white,minimum size=4pt,
		inner sep=0pt]
		
		\draw (0,0) node (x_0) [label=left:$x_0$] {}
		-- ++(225:1.5cm) node (u_0) [label=left:$u_0$] {}
		-- ++(-45:1.5cm) node (x_1) [label=left:$x_1$] {}
		-- (x_0);
		
		\draw (x_1)
		-- ++(225:1.5cm) node (u_1) [label=left:$u_1$] {}
		-- ++(-45:1.5cm) node (x_2) [label=left:$x_2$] {}
		-- (x_1);
        
        \draw (x_2)
		-- ++(225:1.5cm) node (u_2) [label=left:$u_2$] {}
		-- ++(-45:1.5cm) node (x_3) [label=left:$x_3$] {}
		-- (x_2);
        
        \tikzstyle{every node}=[draw,circle,fill=white,minimum size=2pt,
		inner sep=2pt]
		
		\draw (5,-1) node (1) [label=right:{$C_1=\{x_0,u_0,x_1\}$}] {1}
		-- ++(-90:2.2cm) node (2) [label=right:{$C_2=\{x_1,u_1,x_2\}$}] {2}
		-- ++(-90:2.2cm) node (3) [label=right:{$C_3=\{x_2,u_2,x_3\}$}] {3};

		\end{tikzpicture}
		\caption{Sparsity graph for the problem in (\ref{problem:classical_objective}--\ref{problem:classical_dynamics}) to the left and its corresponding clique tree to the
        right.}
		\label{fig:backward-sparsity}

	\end{figure}
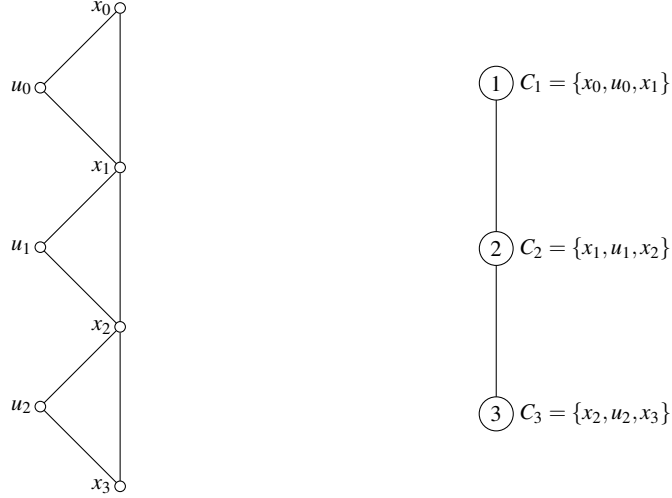

When looking at (\ref{problem:classical_objective}--\ref{problem:classical_dynamics})
it can be put in the almost separable formulation in (\ref{eqPaperIII:CPS}) by defining
\begin{align*}
\bar F_1(x_0, u_0,x_1)&=\mathcal I_{\mathcal D}(x_0)+
\frac{1}{2}\begin{bmatrix}x_0\\u_0\end{bmatrix}^TQ\begin{bmatrix}x_0\\u_0\end{bmatrix}+
\mathcal I_{\mathcal C_0}(x_0, u_0,x_1)\\
\bar F_{k+1}(x_k, u_k,x_{k+1})&=
\frac{1}{2}\begin{bmatrix}x_k\\u_k\end{bmatrix}^TQ\begin{bmatrix}x_k\\u_k\end{bmatrix}+
\mathcal I_{\mathcal C_k}(x_k, u_k,x_{k+1}), 
\quad k=1,\ldots, N-2\\
\bar F_{N}(x_{N-1},u_{N-1},x_N)&=
\frac{1}{2}\begin{bmatrix}x_{N-1}\\u_{N-1}\end{bmatrix}^TQ\begin{bmatrix}x_{N-1}\\u_{N-1}\end{bmatrix}+
\mathcal I_{\mathcal C_{N-1}}(x_{N-1}, u_{N-1},x_N)
+\frac{1}{2}x_N^TSx_N
\end{align*}
where $\mathcal I_{\mathcal C_k}(x_k, u_k,x_{k+1})$ is the indicator function for
the set
$$\mathcal C_k=\left\lbrace (x_k, u_k,x_{k+1})\mid x_{k+1} = Ax_k+Bu_k\right\rbrace$$
and where 
$\mathcal I_{\mathcal D}(x_0)$ is the indicator function for the set
$$\mathcal D=\left\lbrace x_0\mid x_0=\bar x\right\rbrace$$
We have assumed that $v_k=0$. It should be stressed that the derivations done
below easily can be extended to the general case. 
The sparsity graph for this problem is depicted in Figure~\ref{fig:backward-sparsity} for
the case of $N=3$. For ease of 
notation we label the nodes with the states and control signals.\footnote{Here we use a supernode for all 
components of a state and a control signal, respectively. 
In case there is further structure in the dynamic equations
such that not all components of the control signal and the
states are coupled, then more detailed modeling could 
potentially be benificial.} The cliques for this
graph are
$$C_{k+1}=\left\lbrace x_k,u_k,x_{k+1}\right\rbrace,\quad k = 0,\ldots, N-1$$
To obtain a backward dynamic problem formulation we define a clique tree by 
taking $C_1$ as root as seen in Figure~\ref{fig:backward-sparsity}. 
We then assign $\bar F_k$ to $C_k$. This is the only information that has to be provided to a 
general purpose software for solving loosely coupled quadratic programs.

One may of course derive the well-known Riccati-recursion based solution from what has been defined above. The $k$:th problem to solve is 
\begin{align*}
\min_{u_k}\;&\frac{1}{2}\begin{bmatrix}x_k\\u_k\end{bmatrix}^T
Q\begin{bmatrix}x_k\\u_k\end{bmatrix}+m_{k+1,k}(x_{k+1})\\
{\rm s.t.}\;&x_{k+1}=Ax_k+Bu_k
\end{align*}
for given $x_{k}$ starting with $k=N-1$ going down to $k=0$, where $m_{N,N-1}(x_N)=\frac{1}{2}x_N^TSx_N$. 
The optimality conditions are for $k=N-1$
$$
\begin{bmatrix}
S & & I\\
  &Q_2&-B^T\\
  I&-B&
\end{bmatrix}
\begin{bmatrix}
x_N\\u_{N-1}\\\lambda_N
\end{bmatrix}=
\begin{bmatrix}
0\\-Q_{12}^Tx_{N-1}\\Ax_{N-1}
\end{bmatrix}
$$
We notice that no pre-processing of constraints is needed since $\begin{bmatrix}I &-B\end{bmatrix}$ has full row rank. Therefore by the results
in Section 2.5, if the 
overall problem has a unique solution, so does the above problem. 
It is however not easy to give conditions for when the overall problem has
a unique solution. 
The optimality conditions above 
are equivalent to $\lambda_N=-Sx_N$, $x_N=Ax_{N-1}+Bu_{N-1}$ and
the equation
$$G_{N-1}u_{N-1}=-H_{N-1}^Tx_{N-1}$$
where $G_{N-1}=Q_2+B^TSB$ and $H_{N-1}=Q_{12}+A^TSB$. Let 
$F_{N-1}=Q_1+A^TSA$.  Then, if $\begin{bmatrix}A&B\end{bmatrix}$ has full row rank, $Q$ is positive
definite on the nullspace of $\begin{bmatrix}A&B\end{bmatrix}$, and if
$S\succ 0$, it holds that\footnote{Notice that the assumptions are not necessary for the block matrix to be postive definite. Moreover, for the case when it is only positive semidefinite, we still
have a solution $u_{N-1}$, but it is not unique. One may use
pseudo inverse to obtain one solution. This follows from the 
genearlized Schur complement formula. The full row rank assumption 
is equivalent to $(A,B)$ not having any uncontrollable modes corresponding to
zero eigenvalues. The positive definiteness of $Q$ on the nullspace of
$\begin{bmatrix}A&B\end{bmatrix}$ is equivalent to $C(zI-A)^{-1}B+D$ not
having any zeros at the origin where $Q=\begin{bmatrix}C&D\end{bmatrix}^T\begin{bmatrix}C&D\end{bmatrix}$ is a full rank factorization.}  
$$\begin{bmatrix}F_{N-1}&H_{N-1}\\H_{N-1}^T&G_{N-1}\end{bmatrix}=Q+
\begin{bmatrix}A&B\end{bmatrix}^TS\begin{bmatrix}A&B\end{bmatrix}\succ 0$$
Therefore is $G_{N-1}$ positive definite by the Schur complement
formula, and hence there is a unique solution
$$u_{N-1}=-G_{N-1}^{-1}H_{N-1}^Tx_{N-1}$$
Back-substitution of this into the objective function shows that
$$m_{N-1,N-2}(x_{N-1})=\frac{1}{2}x_{N-1}^TS_{N-1}x_{N-1}$$
where $S_{N-1}=F_{N-1}-H_{N-1}G_{N-1}^{-1} H_{N-1}$
which is the first step of the well-known Riccati recursion. Notice
that $S_{N-1}$ is positive definite by the Schur complement formula. 
Repeating the above steps for the remaining problems shows that 
the overall solution can be obtained using the Riccati recursion. 
Notice
that a general purpose solver instead factorizes the local optimality conditions at each step. 
It is well-known that the Riccati recursion provides a factorization of the overall KKT matrix, \cite{van+boy+nou01}. It can be shown that the message 
passing algorithm does the same, \cite{khoshfetrat2016distributed}. 
The main point of this
chapter is, however, that there is no need to derive Riccati
recursions or have any interpretations as factorizations. 
This becomes even more evident when we look at not so 
well-studied MPC formulations, where the corresponding 
structure making it possible to see how Riccati recursions can
be used are only revealed after cumbersome manipulations
of the KKT equations. In some cases it is not even possible
to derive Riccati recursions, which is the point of the 
next subsection. 
\subsection{Forward Dynamic Programming}
Instead of taking $C_1$ as root as we did in the previous subsection we can also 
choose $C_{N}$ as root. We then obtain a forward dynamic programming formulation. We assign the functions to the cliques in the same way. 
The initial problem to solve is
\begin{align*}
\min_{u_0}\;&\frac{1}{2}\begin{bmatrix}x_0\\u_0\end{bmatrix}^TQ
\begin{bmatrix}x_0\\u_0\end{bmatrix}\\
{\rm s.t.}\;&x_1=Ax_0+Bu_0,\quad x_0=\bar x
\end{align*}
parametrically for all possible values of $x_1$. Here we realize that
the constraints for $(x_0,u_0)$ do not satisfy the full row rank 
assumption, i.e. 
$$\begin{bmatrix}I&\\A&B\end{bmatrix}$$
does not have full row rank. Therefore pre-processing is required, 
which can be done using e.g. a QR-factorization on $B$. This will 
result in constraints on $x_1$ that should be passed to the next
problem, and then this procedure should be repeated. Because of this, there is no such clean solution procedure for the forward approach as
for the backward approach, and particularly no Riccati 
recursion based approach. However, the general message passing
approach indeed works.
\subsection{Parallel Computation}
In the previous cases we had a tree that was a chain. It was then 
possible to let either of the end cliques be the root of the tree. 
However, nothing stops us from picking up any one of the midle
cliques as the root. This would result in two branches, and it would then
be possible to solve the problems in the two branches in parallel, one
branch using the backward approach, and one using the forward 
approach. This does however not generalize to more than two parallel
branches. If we want to have three or more we need to proceed 
differently. 

To this end, let us consider a simple example where $N=6$. Let us also assume that we want to solve this problem using two computational agents such that each would perform independently, and hence in parallel. For this, we define dummy variables $\bar{u}_0$ and $\bar{u}_1$ and constrain them as
	
	\begin{equation*}
	\bar{u}_0=x_3, \quad \bar{u}_1=x_6
	\end{equation*}

This is similar to what is done in \cite{nie17} to obtain parallel 
computations. We also define the following sets
	
	\begin{equation}
	\begin{aligned}
	&\mathcal{C}_{-1}=\{x_0:x_0=\bar{x}\} \\
	&\mathcal{C}_k=\{(x_k,u_k,x_{k+1}): \: x_{k+1}=Ax_k+Bu_k\}; \: k=0,1\\
	&\mathcal{C}_2=\{(x_2,u_2,\bar{u}_0): \: \bar{u}_0=Ax_2+Bu_2\}\\
	&\mathcal{C}_k=\{(x_k,u_k,x_{k+1}): \: x_{k+1}=Ax_k+Bu_k\}; \: k=3,4\\
	&\mathcal{C}_5=\{(x_5,u_5,\bar{u}_1): \: \bar{u}_1=Ax_5+Bu_5\}\\
	&\mathcal{D}_0=\{(x_3,\bar{u}_0): \: \bar{u}_0=x_3\}\\
	&\mathcal{D}_1=\{(x_6,\bar{u}_1): \: \bar{u}_1=x_6\}
	\end{aligned}
	\end{equation} 
	
	Then the problem in (\ref{problem:classical_objective}-\ref{problem:classical_dynamics}) can be equivalently written as
	
	\begin{align}
	\underset{u}{\text{min}}
	\:\:&\frac{1}{2}\sum_{k=0}^{1}\begin{bmatrix}x_k\\u_k\end{bmatrix}^T Q \begin{bmatrix}x_k\\u_k\end{bmatrix}
	+\mathcal{I}_{\mathcal{C}_{k}}\{x_k,u_k,x_{k+1}\}+
	\label{example:classical_reconstructed}
	\\
	&\frac{1}{2}\begin{bmatrix}x_2\\u_2\end{bmatrix}^T Q \begin{bmatrix}x_2\\u_2\end{bmatrix}
	+\mathcal{I}_{\mathcal{C}_{2}}\{x_2,u_2,\bar{u}_0\}+
	\nonumber
	\\
	&\frac{1}{2}\sum_{k=3}^{4}\begin{bmatrix}x_k\\u_k\end{bmatrix}^T Q \begin{bmatrix}x_k\\u_k\end{bmatrix}
	+\mathcal{I}_{\mathcal{C}_{k}}\{x_k,u_k,x_{k+1}\}+
	\nonumber
	\\
	&\frac{1}{2}\begin{bmatrix}x_5\\u_5\end{bmatrix}^T Q \begin{bmatrix}x_5\\u_5\end{bmatrix}
	+\mathcal{I}_{\mathcal{C}_{5}}\{x_5,u_5,\bar{u}_1\}+ \frac{1}{2}\bar u_1^TS\bar u_1+
	\nonumber
	\\
	&\mathcal{I}_{\mathcal{C}_{-1}}\{x_0\}+\mathcal{I}_{\mathcal{D}_0}\{x_3,\bar{u}_0\}+\mathcal{I}_{\mathcal{D}_1}\{x_6,\bar{u}_1\}
	\nonumber
	\end{align}
	
\noindent where $\mathcal{I}_{\mathcal{X}}(x)$ is the indicator function
for the set $\mathcal{X}$.
Notice that it is important to define $\mathcal{C}_2$ in terms of $\bar{u}_0$ and not in terms of $x_3$, and similarly for $\mathcal{C}_5$. This trick will allow us to have two independent computational agents. The reason for this will be clear later on.
	
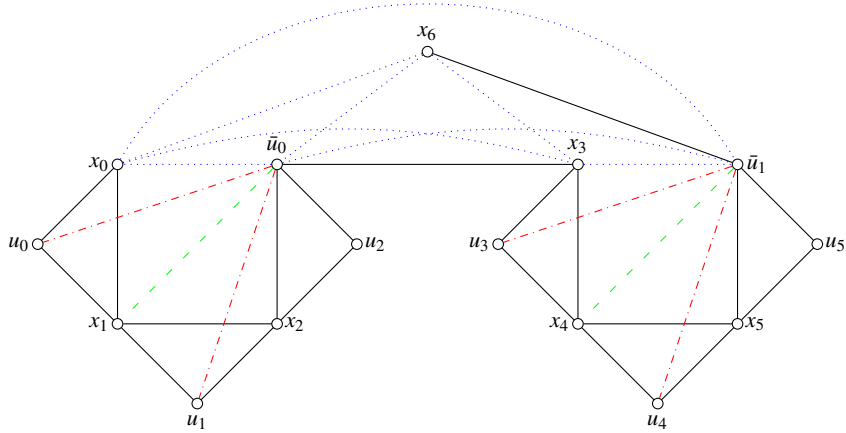
\begin{figure}
		
		\centering
		\begin{tikzpicture}
		\tikzstyle{every node}=[draw,circle,fill=white,minimum size=4pt,
		inner sep=0pt]
		
		\draw (0,0) node (x_0) [label=left:$x_0$] {}
		-- ++(225:1.5cm) node (u_0) [label=left:$u_0$] {}
		-- ++(-45:1.5cm) node (x_1) [label=left:$x_1$] {}
		-- (x_0);
		
		\draw (x_1)
		-- ++(-45:1.5cm) node (u_1) [label=-90:$u_1$] {}
		-- ++(45:1.5cm) node (x_2) [label=right:$x_2$] {}
		-- (x_1);
		
		\draw (x_2)
		-- ++(45:1.5cm) node (u_2) [label=right:$u_2$] {}
		-- ++(135:1.5cm) node (u_0_bar) [label=90:$\bar{u}_0$] {}
		-- (x_2);
		
		\draw (u_0_bar)
		-- ++(0:4cm) node (x_3) [label=90:$x_3$] {}
		;
		
		\draw (x_3)
		-- ++(225:1.5cm) node (u_3) [label=left:$u_3$] {}
		-- ++(-45:1.5cm) node (x_4) [label=left:$x_4$] {}
		-- (x_3);
		
		\draw (x_4)
		-- ++(-45:1.5cm) node (u_4) [label=-90:$u_4$] {}
		-- ++(45:1.5cm) node (x_5) [label=right:$x_5$] {}
		-- (x_4);
		
		\draw (x_5)
		-- ++(45:1.5cm) node (u_5) [label=right:$u_5$] {}
		-- ++(135:1.5cm) node (u_1_bar) [label=right:$\bar{u}_1$] {}
		-- (x_5);
		
		\draw (u_1_bar)
		-- ++(160.0005:4.3858cm) node (x_6) [label=90:$x_6$] {}
		;
		
		
		
		\draw[dotted,blue] (u_0_bar) to (x_0);
		\draw[dotted,blue] (u_0_bar) to [out=15,in=165] (u_1_bar);
		\draw[dotted,blue] (x_3) to (u_1_bar);
		\draw[dotted,blue] (x_0) to [out=15,in=165] (x_3);
		\draw[dotted,blue] (x_0) to [out=60,in=120] (u_1_bar);
		\draw[dotted,blue] (x_6) to (x_3);
		\draw[dotted,blue] (x_6) to (u_0_bar);
		\draw[dotted,blue] (x_6) to (x_0);
		
		\draw[loosely dashed,green] (u_0_bar) to (x_1);
		\draw[loosely dashed,green] (u_1_bar) to (x_4);
		
		\draw[dashdotted,red] (u_0_bar) to (u_0);
		\draw[dashdotted,red] (u_0_bar) to (u_1);
		\draw[dashdotted,red] (u_1_bar) to (u_3);
		\draw[dashdotted,red] (u_1_bar) to (u_4);

		\end{tikzpicture}
		\caption{A modified sparsity graph for the problem in (\ref{example:classical_reconstructed}). The initial sparsity graph, without any modification, is marked with solid lines.}
		\label{SparsityGraphMPCExample}

	\end{figure}

\begin{figure}
		
		\centering
		\begin{tikzpicture}
		
		
		\tikzstyle{every node}=[draw,circle,fill=white,minimum size=2pt,
		inner sep=2pt]
		
		\draw (0,0) node (1) [label=right:{$C_1=\{x_0,\bar{u}_0,x_3,\bar{u}_1,x_6\}$}] {1}
		-- ++(250:1.5cm) node (2) [label=left:{$C_2=\{x_0,u_0,x_1,\bar{u}_0\}$}] {2}
		-- ++(250:2cm) node (4) [label=left:{$C_3=\{x_1,u_1,x_2,\bar{u}_0\}$}] {4}  
		-- ++(250:2cm) node (6) [label=left:{$C_4=\{x_2,u_2,\bar{u}_0\}$}] {6}
		;
		
		\draw (1) 
		-- ++(290:1.5cm) node (3) [label=right:{$C_5=\{x_3,u_3,x_4,\bar{u}_1\}$}] {3}  
		-- ++(290:2cm) node (5) [label=right:{$C_{6}=\{x_4,u_4,x_5,\bar{u}_1\}$}] {5} 
		-- ++(290:2cm) node (7) [label=right:{$C_{7}=\{x_5,u_5,\bar{u}_1\}$}] {7} 
		;

		\end{tikzpicture}
		\caption{Corresponding clique tree for the modified sparsity graph shown in Figure \ref{SparsityGraphMPCExample}.}
		\label{CliqueTreeMPCExample}
		
	\end{figure}
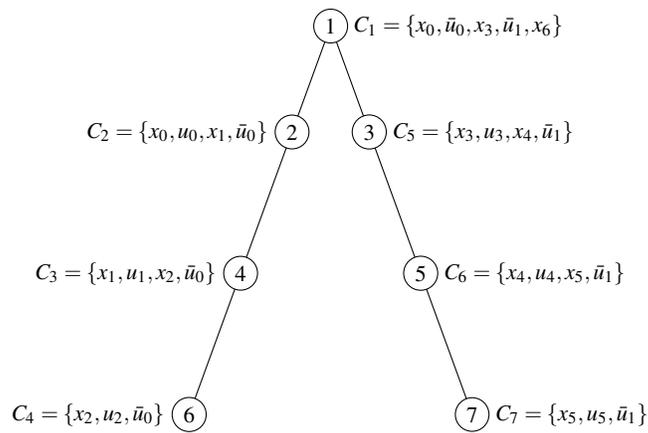

Let us consider the sparsity graph for the problem in (\ref{example:classical_reconstructed}), which is depicted in Figure \ref{SparsityGraphMPCExample}, marked with solid lines. In order to take obtain a clique tree that facilitates parallel computations, we first add edges, marked with dotted lines, between $x_0$, $\bar{u}_0$, $x_3$, $\bar{u}_1$ and $x_6$ such that they form a \emph{maximal} complete subgraph in the graph. The original graph was chordal, but adding the dotted edges destroyed this. Therefore we make a chordal embedding by adding the dashed edges.\footnote{The added edges corresponds to saying that terms 
in the objective function are functions of variables which they are actually not.} We actually add even more edges, which corresponds
to merging cliques. These are the dash-dotted edges. The reason we do this
is that we do not need computational agents for more cliques than the ones
we get after the merging. 
A clique tree which corresponds to the modified sparsity graph in Figure \ref{SparsityGraphMPCExample}, is illustrated in Figure \ref{CliqueTreeMPCExample}. This clique tree obviously enables parallel computations. The different terms in (\ref{example:classical_reconstructed}) are assigned such that rows one and
two are assigned to the left branch, rows three and four to the right branch and the last row to the root. 

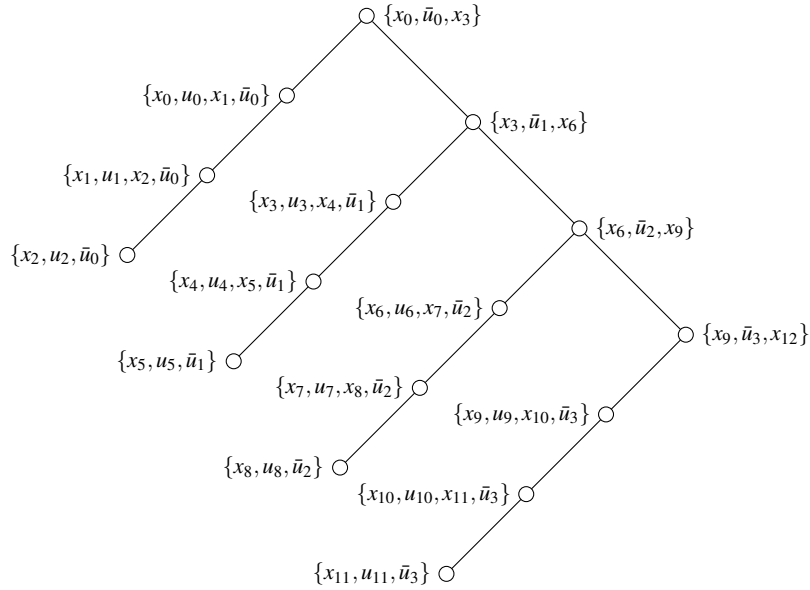
\begin{figure}
		
		\centering
		\begin{tikzpicture}
		
		
		\tikzstyle{every node}=[draw,circle,fill=white,minimum size=2pt,
		inner sep=2pt]
		
		\draw (0,0) node (1) [label=right:{$\{x_0,\bar u_0,x_3\}$}] {}
		-- ++(225:1.5cm) node (2) [label=left:{$\{x_0,u_0,x_1,\bar{u}_0\}$}] {}
		-- ++(225:1.5cm) node (3) [label=left:{$\{x_1,u_1,x_2,\bar{u}_0\}$}] {}  
		-- ++(225:1.5cm) node (4) [label=left:{$\{x_2,u_2,\bar{u}_0\}$}] {}
		;
		
		\draw (1) 
		-- ++(-45:2cm) node (5) [label=right:{$\{x_3,\bar u_1,x_6\}$}] {}  
		-- ++(225:1.5cm) node (6) [label=left:{$\{x_3,u_3,x_4,\bar{u}_1\}$}] {}
		-- ++(225:1.5cm) node (7) [label=left:{$\{x_4,u_4,x_5,\bar{u}_1\}$}] {}  
		-- ++(225:1.5cm) node (8) [label=left:{$\{x_5,u_5,\bar{u}_1\}$}] {}
		;        

		\draw (5) 
		-- ++(-45:2cm) node (9) [label=right:{$\{x_6,\bar u_2,x_9\}$}] {}  
		-- ++(225:1.5cm) node (10) [label=left:{$\{x_6,u_6,x_7,\bar{u}_2\}$}] {}
		-- ++(225:1.5cm) node (11) [label=left:{$\{x_7,u_7,x_8,\bar{u}_2\}$}] {}  
		-- ++(225:1.5cm) node (12) [label=left:{$\{x_8,u_8,\bar{u}_2\}$}] {}
		;        		

		\draw (9) 
		-- ++(-45:2cm) node (13) [label=right:{$\{x_9,\bar u_3,x_{12}\}$}] {}  
		-- ++(225:1.5cm) node (14) [label=left:{$\{x_9,u_9,x_{10},\bar{u}_3\}$}] {}
		-- ++(225:1.5cm) node (15) [label=left:{$\{x_{10},u_{10},x_{11},\bar{u}_3\}$}] {}  
		-- ++(225:1.5cm) node (16) [label=left:{$\{x_{11},u_{11},\bar{u}_3\}$}] {}
		;        		       
        
		\end{tikzpicture}
		\caption{Clique tree with four parallel branches.}
		\label{CliqueTreeMPCExampleExtended}
		
	\end{figure}

Notice that in this particular example we obtained a clique tree with two parallel branches. However, we can generalize to several parallel 
branches by introducing more dummy variables and constraints. Also it is worth pointing out that the subproblem which is assigned to the root of clique tree can be seen as an LQ problem and hence we can use the procedure discussed above recursively. This is 
similar to what is presented in \cite{nie17}. This is obtained in the example 
above by not connecting all of $x_0$, $\bar u_0$, $x_3$, $\bar u_1$ and $x_6$.
Instead one should connect all of $x_0$, $\bar u_0$ and $x_3$ and then 
all of $x_3$, $\bar u_1$ and $x_6$ separately. This will then split the clique $C_1$ 
into two cliques in the clique tree. The dynamics for this LQ problem is 
defined by the sets $\mathcal D_k$, $k=0,1$. The incremental costs for this problem
will be the messages sent by the children of these qliques. Notice that we do not 
really have to know that the resulting problem will be an LQ problem---that is just
an interpretation. We only need to know how to split the root qlicue into one for
each parallel branch. If we want four parallel branches the clique tree will
be like in Figure~\ref{CliqueTreeMPCExampleExtended}.

\subsection{Merging of Cliques}
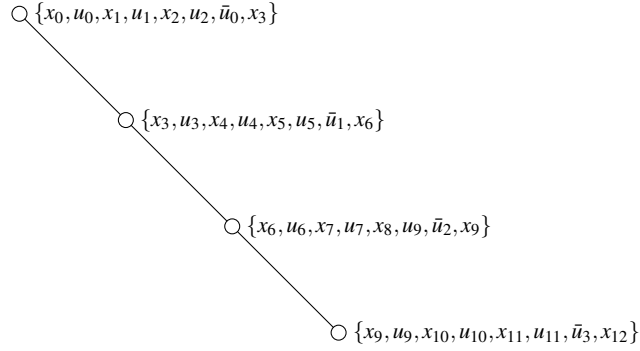
\begin{figure}
		
		\centering
		\begin{tikzpicture}
		
		
		\tikzstyle{every node}=[draw,circle,fill=white,minimum size=2pt,
		inner sep=2pt]
		
		\draw (0,0) node (1) [label=right:{$\{x_0,u_0,x_1,u_1,x_2,u_2,\bar u_0,x_3\}$}] {}
		-- ++(-45:2cm) node (5) [label=right:{$\{x_3,u_3,x_4,u_4,x_5,u_5,\bar u_1,x_6\}$}] {}  
		-- ++(-45:2cm) node (9) [label=right:{$\{x_6,u_6,x_7,u_7,x_8,u_9,\bar u_2,x_9\}$}] {}  
		-- ++(-45:2cm) node (13) [label=right:{$\{x_9,u_9,x_{10},u_{10},x_{11},u_{11},\bar u_3,x_{12}\}$}] {};        		       
        
		\end{tikzpicture}
		\caption{Merged clique tree.}
		\label{CliqueTreeMerged}
		
	\end{figure}

It is not allways the case that one has one agent or processor available for each and
every clique in the clique tree. What then can be done is to merge cliques 
until there are as many cliques as there are processors. Let us consider the
clique tree in Figure~\ref{CliqueTreeMPCExampleExtended}. We can merge the cliques
in each and every parallel branch into one clique. The resulting clique tree will
then be a chain as depicted in Figure~\ref{CliqueTreeMerged}. This could have been
done even before the clique tree was constructed. However, it is benificial for
each of the four agents in the example to utilize the additional structure within
their cliques, i.e. that they have an internal chain structure. This information
would have been lost in case the cliques were merged before the clique tree was
formed. 
\section{Regularized MPC}
A regularized MPC problem is obtained from the above MPC
problem by adding a regularization term to the objective function. 
Typically this is term is either proportional to the squared 
$l^2$ (Euclidian) norm, so-called Tikhonov-regularization, or proportional to
the $l^1$ norm, so-called Lasso-regularization. In both cases convexity
is preserved, since the sum of two convex functions is convex. In the former
case a quadratic objective function will remain quadratic. This is however
not the case for Lasso-regularization. A fairly general Lasso-reguralized 
problem is: 
\begin{align}
\min_{u} &\frac{1}{2}\sum_{k=0}^{N-1}\begin{bmatrix}x_k\\u_k\end{bmatrix}^TQ
\begin{bmatrix}x_k\\u_k\end{bmatrix}+\frac{1}{2}x_N^TSx_N+
\sum_{k=0}^{N-1}\|y_k\|_1\\
{\rm s.t.}\;& x_{k+1}=Ax_k+Bu_k+v_k,\quad x_0=\bar x\\
&Cx_k+Du_k\leq e_k\\
&Ex_k+Fu_k=y_k
\end{align}
\subsection{Equivalent QP}
An equivalent problem formulation for the case of no inequality constraints and $v_k=0$
is:
\begin{align}
\min_{u} &\frac{1}{2}\sum_{k=0}^{N-1}\begin{bmatrix}x_k\\u_k\end{bmatrix}^TQ
\begin{bmatrix}x_k\\u_k\end{bmatrix}+\frac{1}{2}x_N^TSx_N+
\sum_{k=0}^{N-1}t_k\\
{\rm s.t.}\;& x_{k+1}=Ax_k+Bu_k,\quad x_0=\bar x\\
&t_k\geq Ex_k+Fu_k\\
&t_k\geq -Ex_k-Fu_k
\end{align}
This can be put in the almost separable formulation in (\ref{eqPaperIII:CPS}) by defining
\begin{align*}
\bar F_1(x_0, u_0,t_0,x_1)&=\mathcal I_{\mathcal D}(x_0)+
\frac{1}{2}\begin{bmatrix}x_0\\u_0\end{bmatrix}^TQ\begin{bmatrix}x_0\\u_0\end{bmatrix}+t_0+
\mathcal I_{\mathcal C_0}(x_0, u_0,t_0,x_1)\\
\bar F_{k+1}(x_k, u_k,t_k,x_{k+1})&=
\frac{1}{2}\begin{bmatrix}x_k\\u_k\end{bmatrix}^TQ\begin{bmatrix}x_k\\u_k\end{bmatrix}+t_k+
\mathcal I_{\mathcal C_k}(x_k, u_k,t_k,x_{k+1}), 
\quad k\in\naturals_{N-2}\\
\bar F_{N}(x_{N-1},u_{N-1},t_{N-1}x_N)&=
\frac{1}{2}\begin{bmatrix}x_{N-1}\\u_{N-1}\end{bmatrix}^TQ\begin{bmatrix}x_{N-1}\\u_{N-1}\end{bmatrix}+t_{N-1}\\
&+\mathcal I_{\mathcal C_{N-1}}(x_{N-1}, u_{N-1},t_{N-1},x_N)
+\frac{1}{2}x_N^TSx_N
\end{align*}
where $\mathcal I_{\mathcal C_k}(x_k, u_k,t_k,x_{k+1})$ is the indicator function for
the set
$$\mathcal C_k=\left\lbrace (x_k, u_k,x_{k+1})\mid x_{k+1} = Ax_k+Bu_k;\quad t_k\geq Ex_k+Fu_k;\quad 
t_k\geq -Ex_k-Fu_k\right\rbrace$$
and where 
$\mathcal I_{\mathcal D}(x_0)=\left\lbrace x_0\mid x_0=\bar x\right\rbrace$.
It should be stressed that the derivations done
below easily can be extended to the general case. 
The sparsity graph for this problem is very similar to the one for the classical 
formulation. The cliques for this
graph are
$$C_{k+1}=\left\lbrace x_k,u_k,t_k,x_{k+1}\right\rbrace,\quad k = 0,\ldots, N-1$$
To obtain a backward dynamic problem formulation we define a clique tree by 
taking $C_1$ as root similarly as for the classical formulation. We then assign $\bar F_k$ to $C_k$. This is the only information that has to be provided to a 
general purpose software for solving loosely coupled quadratic programs. We can do the forward dynamic programming formulation
as well as a parallel formulation. 
\section{Stochastic MPC}
We will in this section consider a stochastic MPC problem based on a
scenario tree description. Several other authors have investigated how the structure
stemming from scenario trees can be expolitd, e.g. \cite{gon+gro07,mar+15,lei+pot+boc15,fri+17}. The total  
number of scenarios is $M=d^r$, where $d$ is the number of stochastic
events that can take place at each time stage $k$, and where $r$ is the 
number of time stages for which we consider stochastic events to take 
place. The outcome of the stochastic events are the different values of
$A_k^j$, $B_k^j$ and $v_k^j$. Notice that for values of $k<r$ several of 
these quantities are the same. The optimization problem is
\begin{align}
\min_{u}&
\sum_{j=1}^M\omega_j\left(
\frac{1}{2}\sum_{k=0}^{N-1}
\begin{bmatrix}x_k^j\\u_k^j\end{bmatrix}^TQ
\begin{bmatrix}x_k^j\\u_k^j\end{bmatrix}+
\frac{1}{2}(x_N^j)^TSx_N^j\right)\label{obj:stoch}\\
{\rm s.t.}\;& x_{k+1}^j=A_k^jx_k^j+B_k^ju_k^j+v_k^j,\quad x_0^j=\bar x\label{dyn:stoch}\\
&\bar C u = 0\label{ant:stoch}
&
\end{align}
where the index $j$ refers to the $j$:th scenario. Here we define $u=(u^1,u^2,\ldots,u^M)$ with $u^j=(u_0^j,u_1^j,\ldots,u_{N-1}^j)$, and 
$$\bar C=\begin{bmatrix}
C_{1,2}&-C_{1,2}&      &        &        \\
      &C_{2,3}&-C_{2,3}&        &        \\
      &      &\ddots &\ddots &         \\
      &      &       &C_{M-1,M}&-C_{M-1,M}
\end{bmatrix}$$ 
with 
$$C_{j,j+1}=\begin{bmatrix}I&0\end{bmatrix}$$
where $I$ is an identiy matrix of dimenion $m$ times the number of
nodes that scenarios $j$ and $j+1$ have in common. 
The value of $\omega_j$ is the probability of scenario $j$.
The constraint $\bar Cu=0$ is the so-called non-ancipativity constraint. Instead of
saying that each initial state $x_0^j$ is equal to $\bar x$ we instead 
consider the equivalent formulation $x_0^1=\bar x$ and $x_0^j=x_0^{j+1}$, for $j\in\naturals_{M-1}$.

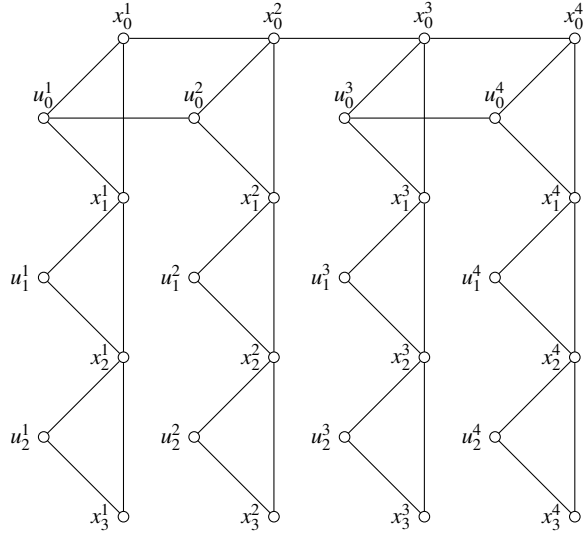
\begin{figure}
		
		\centering
		\begin{tikzpicture}
		\tikzstyle{every node}=[draw,circle,fill=white,minimum size=4pt,
		inner sep=0pt]
		
		\draw (0,0) node (x_01) [label=$x_0^1$] {}
		-- ++(225:1.5cm) node (u_01) [label=$u_0^1$] {}
		-- ++(-45:1.5cm) node (x_11) [label=left:$x_1^1$] {}
		-- (x_01);
		
		\draw (x_11)
		-- ++(225:1.5cm) node (u_11) [label=left:$u_1^1$] {}
		-- ++(-45:1.5cm) node (x_21) [label=left:$x_2^1$] {}
		-- (x_11);
        
        \draw (x_21)
		-- ++(225:1.5cm) node (u_21) [label=left:$u_2^1$] {}
		-- ++(-45:1.5cm) node (x_31) [label=left:$x_3^1$] {}
		-- (x_21);
        
        \draw (2,0) node (x_02) [label=$x_0^2$] {}
		-- ++(225:1.5cm) node (u_02) [label=$u_0^2$] {}
		-- ++(-45:1.5cm) node (x_12) [label=left:$x_1^2$] {}
		-- (x_02);
		
		\draw (x_12)
		-- ++(225:1.5cm) node (u_12) [label=left:$u_1^2$] {}
		-- ++(-45:1.5cm) node (x_22) [label=left:$x_2^2$] {}
		-- (x_12);
        
        \draw (x_22)
		-- ++(225:1.5cm) node (u_22) [label=left:$u_2^2$] {}
		-- ++(-45:1.5cm) node (x_32) [label=left:$x_3^2$] {}
		-- (x_22);        
        
        \draw (4,0) node (x_03) [label=$x_0^3$] {}
		-- ++(225:1.5cm) node (u_03) [label=$u_0^3$] {}
		-- ++(-45:1.5cm) node (x_13) [label=left:$x_1^3$] {}
		-- (x_03);
		
		\draw (x_13)
		-- ++(225:1.5cm) node (u_13) [label=left:$u_1^3$] {}
		-- ++(-45:1.5cm) node (x_23) [label=left:$x_2^3$] {}
		-- (x_13);
        
        \draw (x_23)
		-- ++(225:1.5cm) node (u_23) [label=left:$u_2^3$] {}
		-- ++(-45:1.5cm) node (x_33) [label=left:$x_3^3$] {}
		-- (x_23);        
        
        \draw (6,0) node (x_04) [label=$x_0^4$] {}
		-- ++(225:1.5cm) node (u_04) [label=$u_0^4$] {}
		-- ++(-45:1.5cm) node (x_14) [label=left:$x_1^4$] {}
		-- (x_04);
		
		\draw (x_14)
		-- ++(225:1.5cm) node (u_14) [label=left:$u_1^4$] {}
		-- ++(-45:1.5cm) node (x_24) [label=left:$x_2^4$] {}
		-- (x_14);
        
        \draw (x_24)
		-- ++(225:1.5cm) node (u_24) [label=left:$u_2^4$] {}
		-- ++(-45:1.5cm) node (x_34) [label=left:$x_3^4$] {}
		-- (x_24);        
        
        \draw (x_01) -- (x_02);
        \draw (x_02) -- (x_03);
        \draw (x_03) -- (x_04);
        
        \draw (u_01) -- (u_02);
        \draw (u_03) -- (u_04);
         
		\end{tikzpicture}
		\caption{Sparsity graph for the problem in (\ref{obj:stoch}--\ref{ant:stoch}).}
		\label{fig:stochastic-sparsity}
	\end{figure}

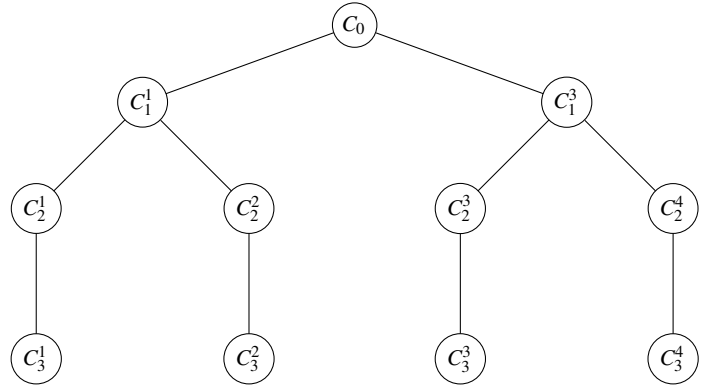
\begin{figure}
		
		\centering
		\begin{tikzpicture}
		
		\tikzstyle{every node}=[draw,circle,fill=white,minimum size=2pt,
		inner sep=2pt]
		
		\draw (0,0) node (0) {$C_0$}
		-- ++(200:3cm) node (11) {$C_1^1$}
		-- ++(225:2cm) node (21) {$C_2^1$}  
		-- ++(270:2cm) node (31) {$C_3^1$}
		;
		
		\draw (0) 
		-- ++(-20:3cm) node (13) {$C_1^3$}  
		-- ++(225:2cm) node (23) {$C_2^3$} 
		-- ++(270:2cm) node (33) {$C_3^3$} 
		;	
        
        \draw (11) 
		-- ++(-45:2cm) node (22) {$C_2^2$}  
		-- ++(270:2cm) node (32) {$C_3^2$} 
		;	
        
        \draw (13) 
		-- ++(-45:2cm) node (24) {$C_2^4$}  
		-- ++(270:2cm) node (34) {$C_3^4$} 
		;	
		\end{tikzpicture}
		\caption{Clique tree for the problem in (\ref{obj:stoch}--\ref{ant:stoch}).}
		\label{fig:stochastic-clique}
		
\end{figure}

We show in Figure~\ref{fig:stochastic-sparsity} the sparsity graph for the case of $d=r=2$ and $N=3$. Then $M=4$. 
We realize that this graph is not chordal. A chordal embedding is obtained by adding
edges such that $C_0=\lbrace x_0^1,x_0^2,x_0^3,x_0^4\rbrace$ is a complete graph. 
Also edges should be added such that $C_1^1=\lbrace x_0^1,u_0^1,x_1^1,x_0^2,u_0^2,x_1^2\rbrace$ and $C_1^3=\lbrace x_0^3,u_0^3,x_1^3,x_0^4,u_0^4,x_1^4\rbrace$ are complete graphs. 
A clique tree for this chordal embedding is shown in Figure~\ref{fig:stochastic-clique},
where $C_{k+1}^j=\lbrace x_k^j,u_k^j,x_{k+1}^j\rbrace$ with $k\in\naturals_{N-1}$
The assignments of functions are for $C_0=\lbrace
x_0^1,x_0^2,x_0^3,x_0^4\rbrace$
$$\bar F_0(x_0^1,x_0^2,x_0^3,x_0^4)=\mathcal I_{\mathcal D}(x_0^1)+
\sum_{j=1}^{M-1}\mathcal I_{\mathcal E}(x_0^j,x_0^{j+1})$$
where $\mathcal D = \lbrace x\mid x=\bar x\rbrace$ and 
$\mathcal E = \lbrace (x,y)\mid x=y\rbrace$.
For $C_1^1$ we assign 
$$\bar F_1^1(x_0^1,u_0^1,x_1^1,x_0^2, u_0^2,x_1^2)=
\sum_{j=1}^2\omega_j\frac{1}{2}\begin{bmatrix}x_0^j\\u_0^j\end{bmatrix}^TQ\begin{bmatrix}x_0^j\\u_0^j\end{bmatrix}+
\mathcal I_{\mathcal C_0^j}(x_0^j, u_0^j,x_1^j)+
+\mathcal I_{\mathcal E}(u_0^1,u_0^2),
$$
for $C_1^3$ we asssign
$$\bar F_1^3(x_0^3,u_0^3,x_1^3,x_0^4, u_0^4,x_1^4)=
\sum_{j=3}^4\omega_j\frac{1}{2}\begin{bmatrix}x_0^j\\u_0^j\end{bmatrix}^TQ\begin{bmatrix}x_0^j\\u_0^j\end{bmatrix}+
\mathcal I_{\mathcal C_0^j}(x_0^j, u_0^j,x_1^j)+
\mathcal I_{\mathcal E}(u_0^3,u_0^4),
$$
for $C_{k+1}^j$, where $k\in\naturals_{N-1}$ and $j\in\naturals_M$, we assign
$$\bar F_{k+1}^j(x_k^j,u_k^j,x_{k+1}^j)=
\omega_j\frac{1}{2}\begin{bmatrix}x_k^j\\u_k^j\end{bmatrix}^TQ\begin{bmatrix}x_k^j\\u_k^j\end{bmatrix}+
\mathcal I_{\mathcal C_k^j}(x_k^j, u_k^j,x_{k+1}^j)$$
and for $C_{N}^j$, where $j\in\naturals_M$, we assign
$$\bar F_{N}^j(x_{N-1}^j,u_{N-1}^j,x_{N}^j)=
\omega_j\frac{1}{2}\begin{bmatrix}x_{N-1}^j\\u_{N-1}^j\end{bmatrix}^TQ\begin{bmatrix}x_{N-1}^j\\u_{N-1}^j\end{bmatrix}+
\omega_j\frac{1}{2}(x_N^j)^TSx_N^j+
\mathcal I_{\mathcal C_{N-1}^j}(x_{N-1}^j, u_{N-1}^j,x_{N}^j)$$
where $\mathcal I_{\mathcal C_k^j}(x_k^j, u_k^j,x_{k+1}^j)$ is the 
indicator function for the set
$$\mathcal C_k^j=\left\lbrace (x_k^j, u_k^j,x_{k+1}^j)\mid x_{k+1}^j = A_k^jx_k^j+B_k^ju_k^j\right\rbrace$$
It is possible to introduce even furhter parallelism by combining the above 
formulation with a parallel formulation in time as described in Section~3.4. 
\section{Distributed MPC}
There are many ways to define distributed MPC problems. We like to think
of them in the following format: 
\begin{align*}
\min_{u} &\frac{1}{2}\sum_{i=1}^M
\sum_{k=0}^{N-1}\begin{bmatrix}x_k^i\\u_k^i\end{bmatrix}^TQ^i
\begin{bmatrix}x_k^i\\u_k^i\end{bmatrix}+\frac{1}{2}\left(x_N^i\right)^T
S^ix_N^i\\
{\rm s.t.}\;& x_{k+1}^i=A^{(i,i)}x_k^i+B^{(i,i)}u_k^i+\sum_{j\in\mathcal N_i}A^{(i,j)}x_k^j+B^{(i,j)}
u_k^j+v_k^i,\quad x_0^i=\bar x^i\\
&C^ix_k^i+D^iu_k^i\leq e_k^i
\end{align*}
for $i\in\naturals_M$, where $\mathcal N_i\subset\naturals_M
\setminus \lbrace i\rbrace$. We see that the only coupling in the problem is in
the dynamic constraints through the summation over $\mathcal N_i$, which 
typically contains few elements. If one consider a sparsity graph for the above
problem one will realize that it is not necessarily chordal. Some heuristic method,
such as presented in \cite{cor:01,kol:09}, can most likely be applied 
sucessfully in many cases to obtain a sparse chordal embedding of the sparsity graph. 
From this a clique tree can be computed using other algorithms presented in \cite{cor:01,kol:09}. See also 
\cite{khoshfetrat2016distributed} for a more detaied
dsicussion on how to compute clique trees. 

We may consider distributed problems
that are stochastic as well. Also extensions to parallelism in time is possible. 
\section{Conclusions}
We have in this chapter shown how it is possible to make use of the inherent 
chordal structure of many MPC formulations in order to exploit IP methods 
that make use of any chordal structure to distribute its computatations over
several computational agents that can work in parallel. We have seen how the 
classical backward Riccait recursion can be seen as a special case of this, albeit
not a parallel recursion, but serial. We have also discussed distributed MPC and
stochastic MPC over scenario trees. The latter formulation can probably be 
extended also to robust MPC over scenario trees. Then the subproblems will
be quadratic feasibility problems and not quadratic programs. Also it should
be possible to consider sum-of-norms regularized MPC. 
We also believe that
it is possible to exploit structure in MPC coming from spatial disretization of PDEs using chordal sparsity. How to carry out these extensions is left for future work. 
\subsection*{Acknowledgements}
The authors want to thank Daniel Axehill and Isak Nielsen for interesting 
discussions regarding parallel computations for Riccati recursions. Shervin 
Parvini Ahmadi has contributed with a figure. 
This reasearch has been supported by WASP, which is gratefully acknowledged. 
%


\bibliographystyle{abbrv}
\bibliography{myref,PaperIII}
\end{document}